\documentclass[12pt]{article}

\usepackage{amsmath, epsfig, cite}
\usepackage{amssymb}
\usepackage{amsfonts}
\usepackage{latexsym}
\usepackage{float}
\usepackage{color}
\usepackage{booktabs}
\usepackage{graphicx}
\usepackage{caption}

\newtheorem{thm}{Theorem}[section]

\newtheorem{cor}[thm]{Corollary}

\newtheorem{lem}[thm]{Lemma}

\numberwithin{equation}{section}

\newcommand{\qed}{{\hfill$\square$}\medskip}
\UseRawInputEncoding
\begin{document}

\begin{center}
{\Large\bf On the binomial transforms of Ap\'ery-like sequences}
\end{center}

\vskip 2mm \centerline{Ji-Cai Liu}
\begin{center}
{\footnotesize Department of Mathematics, Wenzhou University, Wenzhou 325035, PR China\\
{\tt jcliu2016@gmail.com } \\[10pt]
}
\end{center}


\vskip 0.7cm \noindent{\bf Abstract.}
In the proof of the irrationality of $\zeta(3)$ and $\zeta(2)$, Ap\'ery defined two integer sequences through $3$-term recurrences, which are known as the famous Ap\'ery numbers. Zagier, Almkvist--Zudilin and Cooper successively introduced the other $13$ sporadic sequences through variants of Ap\'ery's $3$-term recurrences. All of the $15$ sporadic sequences are called Ap\'ery-like sequences. Motivated by
Gessel's congruences mod $24$ for the Ap\'ery numbers, we investigate the congruences in the form $u_n\equiv \alpha^n \pmod{N_{\alpha}}~(\alpha\in \mathbb{Z},N_{\alpha}\in \mathbb{N}^{+})$ for all of the $15$ Ap\'ery-like sequences $\{u_n\}_{n\ge 0}$. Let $N_{\alpha}$ be the largest positive integer such that $u_n\equiv \alpha^n \pmod{N_{\alpha}}$ for all non-negative integers $n$. We determine the values of $\max\{N_{\alpha}|\alpha \in \mathbb{Z}\}$ for all of the $15$ Ap\'ery-like sequences $\{u_n\}_{n\ge 0}$.
The binomial transforms of Ap\'ery-like sequences provide us a unified approach to this type of congruences for Ap\'ery-like sequences.

\vskip 3mm \noindent {\it Keywords}: Ap\'ery numbers; Ap\'ery-like sequences; binomial transforms; congruences

\vskip 2mm
\noindent{\it MR Subject Classifications}: 11B50, 11B65, 05A10

\section{Introduction}
\subsection{Ap\'ery-like sequences}
In his ingenious proof of the irrationality of $\zeta(3)$ and $\zeta(2)$, Ap\'ery \cite{apery-asterisque-1979} introduced two particular sequences $\{a_n\}_{n\ge 0}$ and $\{b_n\}_{n\ge 0}$ through $3$-term recurrences with given initial values:
\begin{align}
&(n+1)^3 a_{n+1}- (2n+1)(17n^2+17n+5)a_n+n^3 a_{n-1}=0,\quad (a_0=1,a_1=5),\label{a-1}\\[5pt]
&(n+1)^2b_{n+1}-(11n^2+11n+3)b_n-n^2b_{n-1}=0,\quad(b_0=1,b_1=3).\label{a-2}
\end{align}
Ap\'ery also showed that the sequences $\{a_n\}_{n\ge 0}$ and $\{b_n\}_{n\ge 0}$ possess explicit formulas:
\begin{align*}
&a_n=\sum_{k=0}^n{n\choose k}^2{n+k\choose k}^2,\\[5pt]
&b_n=\sum_{k=0}^n{n\choose k}^2{n+k\choose k}.
\end{align*}
These sequences are known as the famous Ap\'ery numbers.

Let $F_{\bf a}(y)$ and $F_{\bf b}(y)$ be the generating functions for the sequences $\{a_n\}_{n\ge 0}$ and $\{b_n\}_{n\ge 0}$, respectively,
\begin{align*}
&F_{\bf a}(y)=\sum_{n=0}^{\infty}a_n y^n,\\
&F_{\bf b}(y)=\sum_{n=0}^{\infty}b_n y^n.
\end{align*}
The recurrences \eqref{a-1} and \eqref{a-2} can be rewritten in forms of differential equations for the generating functions $F_{\bf a}(y)$ and $F_{\bf b}(y)$. In other words, $F_{\bf a}(y)$ is annihilated by the differential operator:
\begin{align*}
\mathcal{L}_{\bf a}=\theta_y^3-y(2\theta_y+1)(17\theta_y^2+17\theta_y+5)+y^2(\theta_y+1)^3,
\end{align*}
and $F_{\bf b}(y)$ is annihilated by the differential operator:
\begin{align*}
\mathcal{L}_{\bf b}=\theta_y^2-y(11\theta_y^2+11\theta_y +3)-y^2(\theta_y+1)^2,
\end{align*}
where $\theta_y=y\frac{d}{dy}$.

It is a surprise that the solutions of the recurrences \eqref{a-1} and \eqref{a-2} are integers.
Motivated by the two interesting examples, Zagier \cite{zagier-b-2009} searched for triples $(A,B,\lambda)\in \mathbb{Z}^3$ such that the solution of the recurrence:
\begin{align}
(n+1)^2u_{n+1}-(An^2+An+\lambda)u_n+Bn^2u_{n-1}=0,\quad (u_{-1}=0,u_0=1)\label{rec-apery-1}
\end{align}
is an integer sequence $\{u_n\}_{n\ge 0}$. The generating function $F(y)=\sum_{n=0}^{\infty}u_n y^n$ for the sequence defined by \eqref{rec-apery-1} is annihilated by the differential operator:
\begin{align}
\mathcal{L}_1=\theta_y^2-y(A\theta_y^2+A\theta_y+\lambda)+By^2(\theta_y+1)^2.\label{diff-ope-1}
\end{align}
In Zagier's computer search, six sporadic solutions are found, the integrality of which are not explained easily. For $(A,B,\lambda)=(11,-1,3)$, we already know that there is a desired solution $\{b_n\}_{n\ge 0}$.

The six sporadic sequences are called Ap\'ery-like sequences of the second kind, which are listed in the following table.
\begin{table}[H]
\caption{\text{Ap\'ery-like sequences of the second kind}}
\centering
\scalebox{0.8}{
\begin{tabular}{cccc}
\toprule
$(A,B,\lambda)$&Name&Other names&Formula\\
\midrule
$(7,-8,2)$&{\bf A}&Franel numbers&$u_n=\sum_{k=0}^n{n\choose k}^3$ \\[7pt]
$(9,27,3)$&{\bf B}&&$u_n=\sum_{k=0}^n(-1)^k3^{n-3k}{n\choose 3k}{3k\choose 2k}{2k\choose k}$ \\[7pt]
$(10,9,3)$&{\bf C}&&$u_n=\sum_{k=0}^n{n\choose k}^2{2k\choose k}$\\[7pt]
$(11,-1,3)$&{\bf D}&Ap\'ery numbers&$u_n=\sum_{k=0}^n{n\choose k}^2{n+k\choose k}$\\[7pt]
$(12,32,4)$&{\bf E}&&$u_n=\sum_{k=0}^n4^{n-2k}{n\choose 2k}{2k\choose k}^2$\\[7pt]
$(17,72,6)$&{\bf F}&&$u_n=\sum_{k=0}^n(-1)^k8^{n-k}{n\choose k}\sum_{j=0}^k{k\choose j}^3$\\
\bottomrule
\end{tabular}
}
\end{table}

Almkvist and Zudilin \cite{az-b-2006} conducted similar searches related to $\{a_n\}_{n\ge 0}$.
They searched for triples $(a,b,c)\in \mathbb{Z}^3$ such that the solution of the recurrence:
\begin{align}
(n+1)^3u_{n+1}-(2n+1)(an^2+an+b)u_n+cn^3u_{n-1}=0,\quad (u_{-1}=0,u_0=1)\label{a-3}
\end{align}
is an integer sequence $\{u_n\}_{n\ge 0}$. Six sporadic sequences $(\delta),(\eta),(\alpha),(\epsilon),(\zeta)$ and $(\gamma)$ were found, which include the
desired sequence $\{a_n\}_{n\ge 0}$.

Cooper \cite{cooper-rj-2012} added a parameter $d$ to the recurrence \eqref{a-3}, and studied a more general recurrence:
\begin{align}
(n+1)^3u_{n+1}-(2n+1)(an^2+an+b)u_n+n(cn^2+d)u_{n-1}=0,\label{rec-apery-2}
\end{align}
with $u_{-1}=0$ and $u_0=1$. Note that the case $d=0$ in \eqref{rec-apery-2} reduces to \eqref{a-3}.
The generating function $F(y)=\sum_{n=0}^{\infty}u_n y^n$ for the sequence defined by \eqref{rec-apery-2} is annihilated by the differential operator:
\begin{align}
\mathcal{L}_2=\theta_y^3-y(2\theta_y+1)(a\theta_y^2+a\theta_y+b)+y^2\left(c(\theta_y+1)^3+d(\theta_y+1)\right).
\label{diff-ope-2}
\end{align}
In Cooper's search, three additional sporadic sequences were found, named $s_7,s_{10}$ and $s_{18}$.

The nine sporadic sequences are called Ap\'ery-like sequences of the first kind, which are listed in the following table.
\begin{table}[H]
\caption{\text{Ap\'ery-like sequences of the first kind}}
\centering
\scalebox{0.8}{
\begin{tabular}{cccc}
\toprule
$(a,b,c,d)$&Name&Other names&Formula\\
\midrule
$(7,3,81,0)$&$(\delta)$&Almkvist-Zudilin numbers&$u_n=\sum_{k=0}^n(-1)^k3^{n-3k}{n\choose 3k}{n+k\choose k}{3k\choose 2k}{2k\choose k}$ \\[7pt]
$(11,5,125,0)$&$(\eta)$&&$u_n=\sum_{k=0}^n(-1)^k{n\choose k}^3({4n-5k-1\choose 3n}+{4n-5k\choose 3n})$ \\[7pt]
$(10,4,64,0)$&$(\alpha)$&Domb numbers&$u_n=\sum_{k=0}^n{n\choose k}^2{2k\choose k}{2n-2k\choose n-k}$ \\[7pt]
$(12,4,16,0)$&$(\epsilon)$&&$u_n=\sum_{k=0}^n{n\choose k}^2{2k\choose n}^2$\\[7pt]
$(9,3,-27,0)$&$(\zeta)$&&$u_n=\sum_{k=0}^n\sum_{l=0}^n{n\choose k}^2{n\choose l}{k\choose l}{k+l\choose n}$ \\[7pt]
$(17,5,1,0)$&$(\gamma)$&Ap\'ery numbers&$u_n=\sum_{k=0}^n{n\choose k}^2{n+k\choose k}^2$ \\[7pt]
$(13,4,-27,3)$&$s_7$&&$u_n=\sum_{k=0}^n{n\choose k}^2{n+k\choose k}{2k\choose n}$ \\[7pt]
$(6,2,-64,4)$&$s_{10}$&Yang-Zudilin numbers&$u_n=\sum_{k=0}^n{n\choose k}^4$ \\[7pt]
$(14,6,192,-12)$&$s_{18}$&&$u_n=\sum_{k=0}^n(-1)^k{n\choose k}{2k\choose k}{2n-2k\choose n-k}({2n-3k-1\choose n}+{2n-3k\choose n})$ \\
\bottomrule
\end{tabular}
}
\end{table}

\subsection{Motivation}
Chowla et al. \cite{ccc-jnt-1980} investigated congruence properties for the Ap\'ery numbers $\{a_n\}_{n\ge 0}$, and conjectured that $a_{2n}\equiv 1\pmod{8},a_{2n+1}\equiv 5\pmod{8},a_{2n}\equiv 1\pmod{3}$ and $a_{2n+1}\equiv 2\pmod{3}$ for all non-negative integers $n$, all of which were proved by
Gessel \cite{gessel-jnt-1982} in the forms:
\begin{align}
&a_n\equiv 5^n\pmod{8},\label{b-1}\\[5pt]
&a_n\equiv (-1)^n \pmod{3}.\label{b-2}
\end{align}
Since $(3,8)=1$, congruences \eqref{b-1} and \eqref{b-2} can be restated in a unified form:
\begin{align}
a_n\equiv 5^n\pmod{24},\label{b-3}
\end{align}
for all non-negative integers $n$.

Another interesting observation is the periodicity for the last digits of the Ap\'ery numbers
$\{b_n\}_{n\ge 0}$. Here are the first few terms of the Ap\'ery numbers $\{b_n\}_{n\ge 0}$:
\begin{align*}
1, 3, 19, 147, 1251, 11253, 104959, 1004307, 9793891, 96918753, 970336269, 9807518757,\cdots
\end{align*}
The last digits of the Ap\'ery numbers $\{b_n\}_{n\ge 0}$ appear to be periodic with $1,3,9,7$, which suggests the conjecture:
\begin{align}
b_n\equiv 3^n\pmod{10},\label{b-4}
\end{align}
for all non-negative integers $n$.

The two interesting examples \eqref{b-3} and \eqref{b-4} motivate us to find more congruences of this type for all of the $15$ Ap\'ery-like sequences $\{u_n\}_{n\ge 0}$.

For an Ap\'ery-like sequence $\{u_n\}_{n\ge 0}$ and an integer $\alpha$, let $N_{\alpha}$ be the largest positive integer such that
\begin{align}
u_n\equiv \alpha^n \pmod{N_{\alpha}},\label{b-5}
\end{align}
for all non-negative integers $n$.

The main purpose of the article is to determine the values of $\beta\in \mathbb{Z}$ such that
$N_{\beta}=\max\{N_{\alpha}|\alpha \in \mathbb{Z}\}$ and the values of $\max\{N_{\alpha}|\alpha \in \mathbb{Z}\}$ for all of the $15$ Ap\'ery-like sequences $\{u_n\}_{n\ge 0}$.

\subsection{Main results}
The binomial transforms of Ap\'ery-like sequences provide us a unified approach to this type of congruence \eqref{b-5} for Ap\'ery-like sequences.
For a sequence $\{u_n\}_{n\ge 0}$ with $u_0=1$ and a complex number $x$, the binomial transform of $\{u_n\}_{n\ge 0}$ is defined by
\begin{align}
v_{n}(x)=\sum_{k=0}^n{n\choose k}(-x)^{n-k} u_k,\label{bino-trans}
\end{align}
for all non-negative integers $n$. Note that $v_{n}(x)$ is a polynomial of degree $n$ with integer coefficients. Throughout the article, given a sequence $\{u_n\}_{n\ge 0}$, the polynomial sequence $\{v_{n}(x)\}_{n\ge 0}$ denotes its corresponding binomial transform \eqref{bino-trans}.

For a positive integer $M$, let $M^{(sf)}$ denote the largest square-free factor of $M$.
For an Ap\'ery-like sequence $\{u_n\}_{n\ge 0}$ and an integer $\alpha$, let
\begin{align}
M_{\alpha}=
\begin{cases}
(v_1(\alpha),v_2(\alpha),v_3(\alpha))\quad&\text{if $\{u_n\}_{n\ge 0}$ is of the second kind},\\[5pt]
(v_1(\alpha),v_2(\alpha),v_3(\alpha),v_4(\alpha))\quad&\text{if $\{u_n\}_{n\ge 0}$ is of the first kind},
\end{cases}\label{def-M}
\end{align}
where $(n_1,n_2,\cdots,n_r)$ denotes the greatest positive common divisor of integers $n_1,n_2,\cdots,n_r$.

The main results of the article consist of the following two theorem.
\begin{thm}\label{t-1}
Let $\{u_n\}_{n\ge 0}$ be one of the $15$ Ap\'ery-like sequences and $\alpha$ be an integer.
For all non-negative integers $n$, we have
\begin{align}
u_n\equiv \alpha^n \pmod{M_{\alpha}^{(sf)}}.\label{main-1}
\end{align}
\end{thm}

\begin{thm}\label{t-2}
Let $\{u_n\}_{n\ge 0}$ be one of the $15$ Ap\'ery-like sequences.
For all integers $\alpha$, $N_{\alpha}$ divides $N_{u_1}$ and
\begin{align}
N_{u_1}=M_{u_1}=
\begin{cases}
(v_2(u_1),v_3(u_1))\quad&\text{if $\{u_n\}_{n\ge 0}$ is of the second kind},\\[5pt]
(v_2(u_1),v_3(u_1),v_4(u_1))\quad&\text{if $\{u_n\}_{n\ge 0}$ is of the first kind},
\end{cases}\label{main-2}
\end{align}
where $N_{\alpha}$ is defined by \eqref{b-5}.
\end{thm}

{\noindent \bf Remark.}
From Theorem \ref{t-2}, we deduce that
\begin{align*}
N_{u_1}=\max\{N_{\alpha}|\alpha \in \mathbb{Z}\}.
\end{align*}
In particular, the congruence $u_n\equiv u_1^n \pmod{N_{u_1}}$
includes all of the congruences $u_n\equiv \alpha^n \pmod{N_{\alpha}}$ as special cases for all integers $\alpha$.

\subsection{Values of $N_{u_1}$}
By \eqref{main-2}, we can determine the values of $N_{u_1}$ for all of the $15$ Ap\'ery-like sequences. Here are the values of $u_1$ and $N_{u_1}$ for the Ap\'ery-like sequences of the second kind.
\begin{table}[H]
\caption{\text{Values of $u_1$ and $N_{u_1}$}}
\centering
\scalebox{0.8}{\label{table-3}
\begin{tabular}{ccccc}
\toprule
Name&Other names&Formula&$u_1$&$N_{u_1}$\\
\midrule
{\bf A}&Franel numbers&$u_n=\sum_{k=0}^n{n\choose k}^3$&$2$&$6$ \\[7pt]
{\bf B}&&$u_n=\sum_{k=0}^n(-1)^k3^{n-3k}{n\choose 3k}{3k\choose 2k}{2k\choose k}$&$3$&$6$ \\[7pt]
{\bf C}&&$u_n=\sum_{k=0}^n{n\choose k}^2{2k\choose k}$&$3$&$6$\\[7pt]
{\bf D}&Ap\'ery numbers&$u_n=\sum_{k=0}^n{n\choose k}^2{n+k\choose k}$&$3$&$10$\\[7pt]
{\bf E}&&$u_n=\sum_{k=0}^n4^{n-2k}{n\choose 2k}{2k\choose k}^2$&$4$&$4$\\[7pt]
{\bf F}&&$u_n=\sum_{k=0}^n(-1)^k8^{n-k}{n\choose k}\sum_{j=0}^k{k\choose j}^3$&$6$&$6$\\
\bottomrule
\end{tabular}
}
\end{table}
From Table \ref{table-3}, we have $b_n\equiv 3^n\pmod{10}$ for all non-negative integers $n$, which
confirms the conjectural congruence \eqref{b-4}.

The values of $u_1$ and $N_{u_1}$ for the Ap\'ery-like sequences of the first kind are listed as follows.
\begin{table}[H]
\caption{\text{Values of $u_1$ and $N_{u_1}$}}
\centering
\scalebox{0.8}{\label{table-4}
\begin{tabular}{ccccc}
\toprule
Name&Other names&Formula&$u_1$&$N_{u_1}$\\
\midrule
$(\delta)$&Almkvist-Zudilin numbers&$u_n=\sum_{k=0}^n(-1)^k3^{n-3k}{n\choose 3k}{n+k\choose k}{3k\choose 2k}{2k\choose k}$&$3$&$24$ \\[7pt]
$(\eta)$&&$u_n=\sum_{k=0}^n(-1)^k{n\choose k}^3({4n-5k-1\choose 3n}+{4n-5k\choose 3n})$&$10$&$10$ \\[7pt]
$(\alpha)$&Domb numbers&$u_n=\sum_{k=0}^n{n\choose k}^2{2k\choose k}{2n-2k\choose n-k}$&$4$&$12$ \\[7pt]
$(\epsilon)$&&$u_n=\sum_{k=0}^n{n\choose k}^2{2k\choose n}^2$&$4$&$24$ \\[7pt]
$(\zeta)$&&$u_n=\sum_{k=0}^n\sum_{l=0}^n{n\choose k}^2{n\choose l}{k\choose l}{k+l\choose n}$&$3$&$6$ \\[7pt]
$(\gamma)$&Ap\'ery numbers&$u_n=\sum_{k=0}^n{n\choose k}^2{n+k\choose k}^2$&$5$&$24$ \\[7pt]
$s_7$&&$u_n=\sum_{k=0}^n{n\choose k}^2{n+k\choose k}{2k\choose n}$&$4$&$8$ \\[7pt]
$s_{10}$&Yang-Zudilin numbers&$u_n=\sum_{k=0}^n{n\choose k}^4$&$2$&$2$ \\[7pt]
$s_{18}$&&$u_n=\sum_{k=0}^n(-1)^k{n\choose k}{2k\choose k}{2n-2k\choose n-k}({2n-3k-1\choose n}+{2n-3k\choose n})$&$12$&$12$ \\
\bottomrule
\end{tabular}
}
\end{table}

\subsection{Structure of the article}
In Section 2, we investigate the binomial transforms of Ap\'ery-like sequences,
which are important ingredients in the proof of Theorems \ref{t-1} and \ref{t-2}. We first establish a key equivalent condition for \eqref{b-5} through binomial transforms of Ap\'ery-like sequences. Secondly, we prove that all of the binomial transforms of Ap\'ery-like sequences satisfy the Gauss congruences through constant term sequences. Thirdly, we establish the generating functions for the binomial transforms of Ap\'ery-like sequences, and the annihilative differential operators for these generating functions. Section 3 is devoted to the proof Theorem \ref{t-1}. In Section 4, we establish an auxiliary result and prove Theorem \ref{t-2} in details.

\section{Binomial transforms}
\subsection{A key equivalent condition}

\begin{lem}\label{equivalent-lem}
Let $\alpha$ be an integer and $N$ be a positive integer. For all non-negative integers $n$,
\begin{align*}
u_n\equiv \alpha^n \pmod{N},
\end{align*}
if and only if for all positive integers $n$,
\begin{align*}
v_n(\alpha)\equiv 0 \pmod{N}.
\end{align*}
\end{lem}
{\noindent\it Proof.}
If $u_n\equiv \alpha^n \pmod{N}$ for all non-negative integers $n$, then
\begin{align*}
v_n(\alpha)\equiv \alpha^{n}\sum_{k=0}^n{n\choose k}(-1)^{n-k}=0\pmod{N},
\end{align*}
for all positive integers $n$.

If $v_n(\alpha)\equiv 0 \pmod{N}$ for all positive integers $n$, then for all non-negative integers $n$,
\begin{align*}
u_n=\sum_{k=0}^n {n\choose k}\alpha^{n-k}v_k(\alpha)\equiv \alpha^n\pmod{N},
\end{align*}
where we have used the binomial inversion formula in the first step.
\qed

\subsection{Gauss congruences}
A sequence $\{u_n\}_{n\ge 1}$ is said to satisfy the Gauss congruences if $u_{np^k}\equiv u_{np^{k-1}}\pmod{p^k}$ for all primes $p$ and all positive integers $n,k$.
A sequence $\{u_n\}_{n\ge 1}$ is called the constant term sequence of a Laurent polynomial $\Lambda\in \mathbb{C}[x_1^{\pm 1},\cdots,x_d^{\pm 1}]$ if $u_n$ is the constant term of $\Lambda^n$ for every $n\ge 1$. It is well-known that Gauss congruences hold for the constant term sequence of a integral Laurent polynomial.

Recently, Gorodetsky \cite{gorodetsky-em-2023} showed that all of the $15$ Ap\'ery-like sequences are constant term sequences.
\begin{lem}(See \cite[Theorem1.1]{gorodetsky-em-2023}.)\label{goro-thm}
Let $\{u_n\}_{n\ge 0}$ be one of the $15$ Ap\'ery-like sequences. Then $\{u_n\}_{n\ge 0}$ is a constant term sequence of a Laurent polynomial with integer coefficients in $2$ or $3$ variables.
\end{lem}

\begin{thm}
Let $\{u_n\}_{n\ge 0}$ be one of the $15$ Ap\'ery-like sequences and $\alpha$ be an integer. Then the corresponding binomial transform $\{v_{n}(\alpha)\}_{n\ge 0}$ satisfies the Gauss congruences:
\begin{align}
v_{np^k}(\alpha)\equiv v_{np^{k-1}}(\alpha)\pmod{p^k},\label{gauss-cong}
\end{align}
for all primes $p$ and all positive integers $n,k$.
\end{thm}
{\noindent \it Proof.}
It suffices to show that the binomial transform $\{v_{n}(\alpha)\}_{n\ge 0}$ is a constant term sequence of a Laurent polynomial with integer coefficients.

For a Laurent series $f(x_1,\cdots,x_k)$, let $CT(f(x_1,\cdots,x_k))$ denote the constant term of the Laurent series. By Lemma \ref{goro-thm}, for all non-negative integers $n$,
\begin{align*}
u_n=CT(\Lambda^n),
\end{align*}
where $\Lambda$ is a Laurent polynomial with integer coefficients in $2$ or $3$ variables.
It follows that for all non-negative integers $n$,
\begin{align*}
&CT\left(\left(\Lambda-\alpha\right)^n\right)\\[5pt]
&=CT\left(\sum_{k=0}^n {n\choose k}(-\alpha)^{n-k}\Lambda^{k}\right)\\[5pt]
&=\sum_{k=0}^n {n\choose k}(-\alpha)^{n-k}CT\left(\Lambda^{k}\right)\\[5pt]
&=\sum_{k=0}^n {n\choose k}(-\alpha)^{n-k}u_k\\[5pt]
&=v_{n}(\alpha).
\end{align*}
We conclude that the binomial transform $\{v_{n}(\alpha)\}_{n\ge 0}$ is a constant term sequence of the integral Laurent polynomial $\Lambda-\alpha$.
\qed

\subsection{Differential operators}
\begin{lem}
For a sequence $\{u_n\}_{n\ge 0}$ and a complex number $x$, let
$F(y)=\sum_{n=0}^{\infty}u_n y^n$ and $G(z)=\sum_{n=0}^{\infty}v_n(x)z^n$.
Then
\begin{align}
G(z)=\frac{1}{1+xz}F\left(\frac{z}{1+xz}\right).\label{gf-1}
\end{align}
\end{lem}
{\noindent\it Proof.}
By \eqref{bino-trans}, we have
\begin{align}
G(z)
&=\sum_{n=0}^{\infty}z^n \sum_{k=0}^n{n\choose k}(-x)^{n-k} u_k\notag\\[7pt]
&=\sum_{k=0}^{\infty}\frac{u_k}{(-x)^k}\sum_{n=k}^{\infty}{n\choose k}(-xz)^n\notag\\[7pt]
&=\sum_{k=0}^{\infty}u_k z^{k}\sum_{n=0}^{\infty}{n+k\choose k}(-xz)^{n}.\label{gf-2}
\end{align}
Note that
\begin{align}
\sum_{n=0}^{\infty}{n+k\choose k}(-xz)^{n}=(1+xz)^{-k-1}.\label{gf-3}
\end{align}
Combining \eqref{gf-2} and \eqref{gf-3}, we complete the proof of \eqref{gf-1}.
\qed

\begin{thm}
For an Ap\'ery-like sequence $\{u_n\}_{n\ge 0}$ and a complex number $x$, let
$F(y)=\sum_{n=0}^{\infty}u_n y^n$ and $G(z)=\sum_{n=0}^{\infty}v_n(x)z^n$.\\
If $F(y)$ is annihilated by the differential operator \eqref{diff-ope-1}, then $G(z)$ is annihilated by the differential operator:
\begin{align}
&\theta_z^2+z\left((3x-A)\theta_z^2+(3x-A)\theta_z+x-\lambda\right)
-z^2(\theta_z+1)^2(2Ax-3x^2-B)\notag\\
&-z^3x(\theta_z+1)(\theta_z+2)(Ax-x^2-B).\label{bino-trans-diff-ope-1}
\end{align}
If $F(y)$ is annihilated by the differential operator \eqref{diff-ope-2}, then $G(z)$ is annihilated by the differential operator:
\begin{align}
&\theta_z^3-z(2\theta_z+1)((a-2x)\theta_z^2+(a-2x)\theta_z+b-x)\notag\\
&-z^2(\theta_z+1)((6ax-6x^2-c)\theta_z^2+(12ax-12x^2-2c)\theta_z+6ax+2bx-7x^2-c-d)\notag\\
&-z^3x(\theta_z+1)(\theta_z+2)(2\theta_z+3)(3ax-2x^2-c)\notag\\
&-z^4x^2(\theta_z+1)(\theta_z+2)(\theta_z+3)(2ax-x^2-c).\label{bino-trans-diff-ope-2}
\end{align}
\end{thm}

{\noindent \it Proof.}
Let
\begin{align}
y=\frac{z}{1+xz}.\label{diff-eq-1}
\end{align}
Then \eqref{gf-1} can be rewritten as
\begin{align}
F(y)=(1+xz)G(z).\label{diff-eq-2}
\end{align}
Differentiating both sides of \eqref{diff-eq-2} with respect to $z$, we obtain
\begin{align}
\theta_yF(y)=\left((1+xz)^2\theta_z+xz(1+xz)\right)G(z).\label{diff-eq-3}
\end{align}
Differentiating both sides of \eqref{diff-eq-3} with respect to $z$, we arrive at
\begin{align}
\theta_y^2F(y)=\left((1+xz)^3\theta_z^2+3xz(1+xz)^2\theta_z+xz(1+xz)(1+2xz)\right)G(z).
\label{diff-eq-4}
\end{align}
Differentiating both sides of \eqref{diff-eq-4} with respect to $z$, we get
\begin{align}
\theta_y^3 F(y)&=((1+xz)^4\theta_z^3+6xz(1+xz)^3\theta_z^2\notag\\
&+xz(1+xz)^2(11xz+4)\theta_z+xz(1+xz)(6(xz)^2+6xz+1))G(z).\label{diff-eq-5}
\end{align}

Finally, substituting \eqref{diff-eq-1}--\eqref{diff-eq-5} into the annihilative differential operators
\eqref{diff-ope-1} and \eqref{diff-ope-2} for $F(y)$, after simplifications, we are led to the annihilative differential operators \eqref{bino-trans-diff-ope-1} and \eqref{bino-trans-diff-ope-2} for $G(z)$.
\qed

From the annihilative differential operators \eqref{bino-trans-diff-ope-1} and \eqref{bino-trans-diff-ope-2} for $G(z)$, we immediately get the recurrences for
the binomial transforms $\{v_{n}(\alpha)\}_{n\ge 0}$ of Ap\'ery-like sequences $\{u_n\}_{n\ge 0}$.
\begin{cor}
Let $\{u_n\}_{n\ge 0}$ be an Ap\'ery-like sequence and $\alpha$ be an integer.\\
If $\{u_n\}_{n\ge 0}$ is
of the second kind, then $\{v_n(\alpha)\}_{n\ge 0}$ satisfies a recurrence in the form:
\begin{align}
n^2v_n(\alpha)+f_1(n)v_{n-1}(\alpha)+f_2(n)v_{n-2}(\alpha)+f_3(n)v_{n-3}(\alpha)=0,\label{bino-trans-rec-1}
\end{align}
where $f_i(x)$ is a polynomial of degree $2$ with integer coefficients for $i=1,2,3$.\\
If $\{u_n\}_{n\ge 0}$ is of the first kind, then $\{v_n(\alpha)\}_{n\ge 0}$ satisfies a recurrence in the form:
\begin{align}
n^3v_n(\alpha)+g_1(n)v_{n-1}(\alpha)+g_2(n)v_{n-2}(\alpha)+
g_3(n)v_{n-3}(\alpha)+g_4(n)v_{n-4}(\alpha)=0,\label{bino-trans-rec-2}
\end{align}
where $g_i(x)$ is a polynomial of degree $3$ with integer coefficients for $i=1,2,3,4$.
\end{cor}

\section{Proof of Theorem \ref{t-1}}
We shall prove the case that $\{u_n\}_{n\ge 0}$ is an Ap\'ery-like sequence of the second kind. The proof for the Ap\'ery-like sequences of the first kind runs analogously, and we omit the details.

Assume that $\{u_n\}_{n\ge 0}$ is an Ap\'ery-like sequence of the second kind. We shall prove that
for all non-negative integers $n$,
\begin{align*}
u_n\equiv \alpha^n \pmod{(v_1(\alpha),v_2(\alpha),v_3(\alpha))^{(sf)}},
\end{align*}
by Lemma \ref{equivalent-lem}, which is equivalent to
\begin{align}
v_n(\alpha)\equiv 0 \pmod{(v_1(\alpha),v_2(\alpha),v_3(\alpha))^{(sf)}},\label{proof-thm1-1}
\end{align}
for all positive integers $n$.

Since $(v_1(\alpha),v_2(\alpha),v_3(\alpha))^{(sf)}$ is square-free, let
$(v_1(\alpha),v_2(\alpha),v_3(\alpha))^{(sf)}=p_1 p_2\cdots p_s$ where $p_1,p_2,\cdots,p_s$ are distinct primes.

It is clear that \eqref{proof-thm1-1} holds for $n=1,2,3$. Let $i\in \{1,2,\cdots,s\}$.
If $p_i$ divides $n$, then, by \eqref{gauss-cong}, we have
\begin{align}
v_n(\alpha)\equiv v_{n/p_i}(\alpha)\pmod{p_i}.\label{proof-thm1-2}
\end{align}
If $p_i$ does not divide $n$, then $1/n^2$ is a $p_i$-adic integer.
It follows from \eqref{bino-trans-rec-1} that
\begin{align}
v_n(\alpha)\equiv c_1v_{n-1}(\alpha)+c_2v_{n-2}(\alpha)+c_3v_{n-3}(\alpha)\pmod{p_i},\label{proof-thm1-3}
\end{align}
where $c_1,c_2,c_3$ are $p_i$-adic integers.

Using \eqref{proof-thm1-2}, \eqref{proof-thm1-3} and the induction on $n$, we obtain
\begin{align*}
v_n(\alpha)\equiv 0\pmod{p_i},
\end{align*}
for all positive integers $n$.

Since $p_1,p_2,\cdots,p_s$ are distinct primes, we have
\begin{align*}
v_n(\alpha)\equiv 0\pmod{p_1p_2\cdots p_s},
\end{align*}
for all positive integers $n$. This completes the proof of \eqref{proof-thm1-1}.

\section{Proof of Theorem \ref{t-2}}
\subsection{An auxiliary result}
\begin{lem}
Let $c$ be an integer and $f(x)$ be a polynomial with integer coefficients. For all integers $\alpha$, we have
\begin{align}
\left(-\alpha+c,f(\alpha)\right)=\left(-\alpha+c,f(c)\right).\label{auxiliary-2}
\end{align}
\end{lem}
{\noindent \it Proof.}
Assume that $f(x)=a_s x^s+a_{s-1}x^{s-1}+\cdots+a_1 x+a_0\in \mathbb{Z}[x]$.
Note that
\begin{align}
f(x)-f(c)=a_s (x^s-c^s)+a_{s-1}(x^{s-1}-c^{s-1})+\cdots+a_1(x-c).\label{auxiliary-1}
\end{align}
We rewrite \eqref{auxiliary-1} in the form:
\begin{align*}
f(x)=(-x+c)g(x)+f(c),
\end{align*}
where $g(x)$ is a polynomial with integer coefficients. It follows that
\begin{align*}
\left(-\alpha+c,f(\alpha)\right)=\left(-\alpha+c,(-\alpha+c)g(\alpha)+f(c)\right)=\left(-\alpha+c,f(c)\right),
\end{align*}
as desired.
\qed

\subsection{Proof of Theorem \ref{t-2}}
Assume that $\{u_n\}_{n\ge 0}$ is an Ap\'ery-like sequence of the second kind.
Using the fact $v_1(\alpha)=-\alpha+u_1$ and \eqref{auxiliary-2}, we have
\begin{align*}
M_{\alpha}&=(v_1(\alpha),v_2(\alpha),v_3(\alpha))\\[5pt]
&=\left((v_1(\alpha),v_2(\alpha)),(v_1(\alpha),v_3(\alpha))\right)\\[5pt]
&=\left((-\alpha+u_1,v_2(\alpha)),(-\alpha+u_1,v_3(\alpha))\right)\\[5pt]
&=\left((-\alpha+u_1,v_2(u_1)),(-\alpha+u_1,v_3(u_1))\right)\\[5pt]
&=\left(-\alpha+u_1,v_2(u_1),v_3(u_1)\right)\\[5pt]
&=\left(-\alpha+u_1,(v_2(u_1),v_3(u_1))\right),
\end{align*}
and so $M_{u_1}=(v_2(u_1),v_3(u_1))$ and $M_{\alpha}$ divides $M_{u_1}$ for all integers $\alpha$.

Similarly, we can show that for an Ap\'ery-like sequence $\{u_n\}_{n\ge 0}$ of the first kind,
\begin{align*}
M_{\alpha}=\left(-\alpha+u_1,(v_2(u_1),v_3(u_1),v_4(u_1))\right),
\end{align*}
and so $M_{u_1}=(v_2(u_1),v_3(u_1),v_4(u_1))$ and $M_{\alpha}$ divides $M_{u_1}$ for all integers $\alpha$.
It follows that for all integers $\alpha$,
\begin{align}
M_{\alpha}{\big |}M_{u_1},\label{proof-thm2-1}
\end{align}
and
\begin{align*}
M_{u_1}=
\begin{cases}
(v_2(u_1),v_3(u_1))\quad&\text{if $\{u_n\}_{n\ge 0}$ is of the second kind},\\[5pt]
(v_2(u_1),v_3(u_1),v_4(u_1))\quad&\text{if $\{u_n\}_{n\ge 0}$ is of the first kind}.
\end{cases}
\end{align*}

By Lemma \ref{equivalent-lem}, we conclude that $N_{\alpha}$ is the largest positive integer such that
\begin{align}
v_n(\alpha)\equiv 0 \pmod{N_{\alpha}},\label{proof-thm2-2}
\end{align}
for all positive integers $n$. By \eqref{proof-thm2-2} and the definition \eqref{def-M}, we have
\begin{align}
N_{\alpha}{\big |} M_{\alpha},\label{proof-thm2-3}
\end{align}
for all integers $\alpha$. Combining \eqref{proof-thm2-1} and \eqref{proof-thm2-3}, we obtian
\begin{align}
N_{\alpha}{\big |} M_{u_1},\label{proof-thm2-4}
\end{align}
for all integers $\alpha$. It remains to prove that
\begin{align}
N_{u_1}=M_{u_1}.\label{proof-thm2-5}
\end{align}

By \eqref{proof-thm2-4}, we have $N_{u_1} |M_{u_1}$. To prove \eqref{proof-thm2-5}, it suffices to show that $M_{u_1}| N_{u_1}$, which is equivalent to
\begin{align}
u_n\equiv u_1^n \pmod{M_{u_1}},\label{proof-thm2-6}
\end{align}
for all non-negative integers $n$. Next, we shall prove that \eqref{proof-thm2-6} is true for all of the $15$ Ap\'ery-like sequences.

By \eqref{main-1}, we have
\begin{align}
u_n\equiv u_1^n \pmod{M_{u_1}^{(sf)}},\label{proof-thm2-7}
\end{align}
for all non-negative integers $n$. If $M_{u_1}$ is square-free, then $M_{u_1}=M_{u_1}^{(sf)}$, and so
\eqref{proof-thm2-6} obviously holds. It follows from the values of $M_{u_1}$ in Tables \ref{table-3} and \ref{table-4} that \eqref{proof-thm2-6} is true for Ap\'ery-like sequences {\bf A}, {\bf B}, {\bf C}, {\bf D}, {\bf F}, ($\eta$), ($\zeta$) and $s_{10}$.

Next, we shall prove that \eqref{proof-thm2-6} is also true for the remaining seven Ap\'ery-like sequences.
Let $M_{u_1}^{*}$ denote the square-free part in the prime factorization of $M_{u_1}$. Since $M_{u_1}^{*}$ divides $M_{u_1}^{(sf)}$, by \eqref{proof-thm2-7}, we have
\begin{align}
u_n\equiv u_1^n \pmod{M_{u_1}^{*}},\label{proof-thm2-8}
\end{align}
for all non-negative integers $n$. Moreover, we have
\begin{align}
(M_{u_1}^{*},M_{u_1}/M_{u_1}^{*})=1.\label{proof-thm2-9}
\end{align}
To prove \eqref{proof-thm2-6}, by \eqref{proof-thm2-8} and \eqref{proof-thm2-9} it suffices to show that
\begin{align}
u_n\equiv u_1^n \pmod{M_{u_1}/M_{u_1}^{*}},\label{proof-thm2-10}
\end{align}
for all non-negative integers $n$.

For the remaining seven Ap\'ery-like sequences, \eqref{proof-thm2-10} reads
\begin{align}
&\sum_{k=0}^n{n\choose k}^2{n+k\choose k}^2\equiv 5^n\pmod{8},\label{cong-1}\\[7pt]
&\sum_{k=0}^n4^{n-2k}{n\choose 2k}{2k\choose k}^2\equiv 4^n \pmod{4},\label{cong-2}\\[7pt]
&\sum_{k=0}^n{n\choose k}^2{2k\choose k}{2n-2k\choose n-k}\equiv 4^n\pmod{4},\label{cong-3}\\[7pt]
&\sum_{k=0}^n(-1)^k{n\choose k}{2k\choose k}{2n-2k\choose n-k}\left({2n-3k-1\choose n}+{2n-3k\choose n}\right)\equiv 12^n \pmod{4},\label{cong-4}\\[7pt]
&\sum_{k=0}^n{n\choose k}^2{2k\choose n}^2\equiv 4^n\pmod{8},\label{cong-5}\\[7pt]
&\sum_{k=0}^n{n\choose k}^2{n+k\choose k}{2k\choose n}\equiv 4^n \pmod{8},\label{cong-6}\\[7pt]
&\sum_{k=0}^n(-1)^k3^{n-3k}{n\choose 3k}{n+k\choose k}{3k\choose 2k}{2k\choose k}
\equiv 3^n\pmod{8}.\label{cong-7}
\end{align}

We remark that Gessel \cite[Theorem 3 (ii)]{gessel-jnt-1982} proved \eqref{cong-1}. Next, we shall prove \eqref{cong-2}--\eqref{cong-7} respectively.

{\noindent \it Proof of \eqref{cong-2}.}
It suffices to show that for all positive integers $n$,
\begin{align}
\sum_{k=0}^n4^{n-2k}{n\choose 2k}{2k\choose k}^2\equiv 0 \pmod{4}.\label{cong-proof-1}
\end{align}
If $n$ is an odd positive integer, then \eqref{cong-proof-1} clearly holds.
If $n$ is an even positive integer, then
\begin{align*}
\sum_{k=0}^{n}4^{n-2k}{n\choose 2k}{2k\choose k}^2\equiv {n\choose n/2}^2=4{n-1\choose n/2-1}^2\equiv 0 \pmod{4},
\end{align*}
as desired.
\qed

{\noindent \it Proof of \eqref{cong-3}.}
It suffices to show that for all positive integers $n$,
\begin{align*}
\sum_{k=0}^n{n\choose k}^2{2k\choose k}{2n-2k\choose n-k}\equiv 0\pmod{4}.
\end{align*}
Since ${2k\choose k}=2{2k-1\choose k-1}\equiv 0\pmod{2}$ for all positive integers $k$, we have
\begin{align*}
\sum_{k=0}^n{n\choose k}^2{2k\choose k}{2n-2k\choose n-k}\equiv 2{2n\choose n}\equiv 0\pmod{4},
\end{align*}
as desired.
\qed

{\noindent \it Proof of \eqref{cong-4}.}
It suffices to show that for all positive integers $n$,
\begin{align*}
\sum_{k=0}^n(-1)^k{n\choose k}{2k\choose k}{2n-2k\choose n-k}\left({2n-3k-1\choose n}+{2n-3k\choose n}\right)\equiv 0 \pmod{4}.
\end{align*}
Since ${2k\choose k}\equiv 0\pmod{2}$ for all positive integers $k$, we have
\begin{align*}
&\sum_{k=0}^n(-1)^k{n\choose k}{2k\choose k}{2n-2k\choose n-k}\left({2n-3k-1\choose n}+{2n-3k\choose n}\right)\\[5pt]
&\equiv {2n\choose n}\left({2n-1\choose n}+{2n\choose n}\right)+(-1)^n{2n\choose n}\left({-n-1\choose n}+{-n\choose n}\right)\\[5pt]
&=2{2n\choose n}\left({2n-1\choose n}+{2n\choose n}\right)\\[5pt]
&\equiv 0 \pmod{4},
\end{align*}
as desired.
\qed

{\noindent \it Proof of \eqref{cong-5}.}
It is easy to check that \eqref{cong-5} holds for $n=0,1,2$.
Let
\begin{align*}
a_n=\sum_{k=0}^n{n\choose k}^2{2k\choose n}^2.
\end{align*}
It suffices to show that for all integers $n\ge 3$,
\begin{align}
a_n\equiv 0\pmod{8}.\label{cong-proof-new-1}
\end{align}

By the recurrence \eqref{rec-apery-2}, we have
\begin{align*}
(2n+1)^3a_{2n+1}-4(4n+1)(12n^2+6n+1)a_{2n}+128n^3a_{2n-1}=0.
\end{align*}
It follows that
\begin{align*}
(4n^2+6n+1)a_{2n+1}+4a_{2n}\equiv 0\pmod{8},
\end{align*}
and so
\begin{align}
a_{2n+1}\equiv \frac{4}{4n^2+6n+1}a_{2n}\equiv 4a_{2n}\pmod{8}.\label{cong-proof-2}
\end{align}

On the other hand, we have
\begin{align}
a_{2n}&=\sum_{k=0}^{2n}{2n\choose k}^2{2k\choose 2n}^2\notag\\[5pt]
&=\sum_{k=0}^{n}{2n\choose 2k}^2{4k\choose 2n}^2+\sum_{k=0}^{n}{2n\choose 2k+1}^2{2(2k+1)\choose 2n}^2\notag\\[5pt]
&=\sum_{k=0}^{n}{2n\choose 2k}^2{4k\choose 2n}^2+4n^2\sum_{k=0}^{n}\frac{1}{(2k+1)^2}{2n-1\choose 2k}^2{2(2k+1)\choose 2n}^2\notag\\[5pt]
&\equiv \sum_{k=0}^{n}{n\choose k}^2{2k\choose n}^2+4n^2\sum_{k=0}^{n}{2n-1\choose 2k}^2{2k+1\choose n}^2\pmod{8},\label{cong-proof-3}
\end{align}
where we have used the following congruence (see \cite{kaz-bsmg-1968}):
\begin{align}
{2s\choose 2t}\equiv {s\choose t} \pmod{4}.\label{cong-proof-4}
\end{align}
Furthermore, for $n\ge 2$ we have
\begin{align}
&n^2\sum_{k=0}^{n}{2n-1\choose 2k}^2{2k+1\choose n}^2\notag\\[5pt]
&=\sum_{k=0}^{n}(2k+1)^2{2(n-1)+1\choose 2k}^2{2k\choose n-1}^2\notag\\[5pt]
&\equiv \sum_{k=0}^{n-1}{n-1\choose k}^2{2k\choose n-1}^2\notag\\[5pt]
&= 4\sum_{k=1}^{n-1}{n-2\choose k-1}^2{2k-1\choose n-2}^2\notag\\[5pt]
&\equiv 0\pmod{2},\label{cong-proof-5}
\end{align}
where we have used \eqref{cong-proof-4}. Finally, combining \eqref{cong-proof-3} and \eqref{cong-proof-5}, we arrive at
\begin{align}
a_{2n}\equiv a_n\pmod{8},\label{cong-proof-6}
\end{align}
for all integers $n\ge 2$.

By using the induction on $n$, \eqref{cong-proof-2} and \eqref{cong-proof-6}, we complete the proof of \eqref{cong-proof-new-1}.
\qed

{\noindent \it Proof of \eqref{cong-6}.}
It is easy to check that \eqref{cong-6} holds for $n=0,1,2,3$.
Let
\begin{align*}
b_n=\sum_{k=0}^n{n\choose k}^2{n+k\choose k}{2k\choose n}.
\end{align*}
It suffices to show that for all integers $n\ge 4$,
\begin{align}
b_n\equiv 0\pmod{8}.\label{cong-proof-7}
\end{align}

By the recurrence \eqref{rec-apery-2}, we have
\begin{align*}
(2n+1)^3b_{2n+1}-2(4n+1)(26n^2+13n+2)b_{2n}-6n(6n-1)(6n+1)b_{2n-1}=0.
\end{align*}
It follows that
\begin{align*}
(4n^2+6n+1)b_{2n+1}+2(2n^2+3n+2)b_{2n}+6nb_{2n-1}\equiv 0\pmod{8},
\end{align*}
which can be rewritten in the form:
\begin{align}
b_{2n+1}\equiv \beta_1 b_{2n}+\beta_2 b_{2n-1}\pmod{8},\label{cong-proof-8}
\end{align}
where $\beta_1,\beta_2$ are $2$-adic integers.

On the other hand, we have
\begin{align}
b_{2n}&=\sum_{k=0}^{2n}{2n\choose k}^2{2n+k\choose k}{2k\choose 2n}\notag\\[5pt]
&=\sum_{k=0}^{n}{2n\choose 2k}^2{2n+2k\choose 2k}{4k\choose 2n}
+\sum_{k=0}^{n}{2n\choose 2k+1}^2{2n+2k+1\choose 2n}{2(2k+1)\choose 2n}\notag\\[5pt]
&=\sum_{k=0}^{n}{2n\choose 2k}^2{2n+2k\choose 2k}{4k\choose 2n}
+4n^2\sum_{k=0}^{n}\frac{1}{(2k+1)^2}{2n-1\choose 2k}^2{2n+2k+1\choose 2n}{2(2k+1)\choose 2n}\notag\\[5pt]
&\equiv \sum_{k=0}^{n}{n\choose k}^2{2n+2k\choose 2k}{4k\choose 2n}
+4n^2\sum_{k=0}^{n}{2n-1\choose 2k}^2{2n+2k+1\choose 2n}{2(2k+1)\choose 2n}\pmod{8},
\label{cong-proof-9}
\end{align}
where we have used \eqref{cong-proof-4}.

Since for $n\ge 1$ and $k\ge 1$,
\begin{align*}
{n\choose k}{4k\choose 2n}=2{n-1\choose k-1}{4k-1\choose 2n-1}\equiv 0\pmod{2},
\end{align*}
and
\begin{align*}
{n\choose k}{n+k\choose k}={n+k\choose 2k}{2k\choose k}=2{n+k\choose 2k}{2k-1\choose k-1}\equiv 0\pmod{2},
\end{align*}
we have
\begin{align}
&\sum_{k=0}^{n}{n\choose k}^2{2n+2k\choose 2k}{4k\choose 2n}\notag\\[5pt]
&\equiv \sum_{k=0}^{n}{n\choose k}^2{n+k\choose k}{4k\choose 2n}\notag\\[5pt]
&\equiv \sum_{k=0}^{n}{n\choose k}^2{n+k\choose k}{2k\choose n}\pmod{8},
\label{cong-proof-10}
\end{align}
where we have also used \eqref{cong-proof-4}.

Furthermore, for $n\ge 2$ we have
\begin{align}
&n^2\sum_{k=0}^{n}{2n-1\choose 2k}^2{2n+2k+1\choose 2n}{2(2k+1)\choose 2n}\notag\\[5pt]
&=n^2\sum_{k=0}^{n}{2(n-1)+1\choose 2k}^2{2n+2k+1\choose 2n}{2(2k+1)\choose 2n}\notag\\[5pt]
&\equiv n^2\sum_{k=0}^{n-1}{n-1\choose k}^2{n+k\choose n}{2k+1\choose n}\notag\\[5pt]
&=2n\sum_{k=0}^{n-1}(2k+1){n-1\choose k}{n+k\choose n}{n-2\choose k-1}{2k-1\choose n-2}\notag\\[5pt]
&\equiv 0\pmod{2}.\label{cong-proof-11}
\end{align}

Finally, combining \eqref{cong-proof-9}--\eqref{cong-proof-11} gives
\begin{align}
b_{2n}\equiv b_n\pmod{8},\label{cong-proof-12}
\end{align}
for all integers $n\ge 2$.

By using the induction on $n$, \eqref{cong-proof-8} and \eqref{cong-proof-12}, we complete the proof of \eqref{cong-proof-7}.
\qed

{\noindent \it Proof of \eqref{cong-7}.}
Let
\begin{align*}
c_n=\sum_{k=0}^n(-1)^k3^{n-3k}{n\choose 3k}{n+k\choose k}{3k\choose 2k}{2k\choose k}.
\end{align*}
Since ${2k\choose k}\equiv 0\pmod{2}$ for $k\ge 1$, we have
\begin{align}
c_n\equiv 3^{n}\sum_{k=0}^n{n\choose 3k}{n+k\choose k}{3k\choose 2k}{2k\choose k}\pmod{8}.
\label{cong-proof-13}
\end{align}

By the recurrence \eqref{rec-apery-2}, we have
\begin{align*}
(2n+1)^3c_{2n+1}-(4n+1)(28n^2+14n+3)c_{2n}+648n^3c_{2n-1}=0.
\end{align*}
It follows that
\begin{align*}
(4n^2+6n+1)c_{2n+1}+(4n^2+6n+5)c_{2n}\equiv 0\pmod{8},
\end{align*}
and so
\begin{align}
c_{2n+1}&\equiv -\frac{4n^2+6n+5}{4n^2+6n+1}c_{2n}\notag\\[5pt]
&=-\left(1+\frac{4}{4n^2+6n+1}\right)c_{2n}\notag\\[5pt]
&\equiv 3c_{2n}\pmod{8}.\label{cong-proof-14}
\end{align}

On the other hand, by \eqref{cong-proof-13} we have
\begin{align}
c_{2n}&\equiv 3^{2n}\sum_{k=0}^{2n}{2n\choose 3k}{2n+k\choose k}{3k\choose 2k}{2k\choose k}\notag\\[5pt]
&=3^{2n}\sum_{k=0}^{n}{2n\choose 6k}{2n+2k\choose 2k}{6k\choose 4k}{4k\choose 2k}\notag\\[5pt]
&+3^{2n}\sum_{k=0}^{n}{2n\choose 3(2k+1)}{2n+2k+1\choose 2k+1}{3(2k+1)\choose 2(2k+1)}{2(2k+1)\choose 2k+1}\pmod{8}.\label{cong-proof-15}
\end{align}

Since ${4k\choose 2k}\equiv 0\pmod{2}$ and ${n\choose 3k}{n+k\choose k}{3k\choose 2k}={2k\choose k}{4k\choose 2k}{n+k\choose 4k}\equiv 0\pmod{4}$ for $k\ge 1$, we have
\begin{align}
&3^{2n}\sum_{k=0}^{n}{2n\choose 6k}{2n+2k\choose 2k}{6k\choose 4k}{4k\choose 2k}\notag\\[5pt]
&\equiv 3^{2n}\sum_{k=0}^{n}{n\choose 3k}{n+k\choose k}{3k\choose 2k}{4k\choose 2k}\notag\\[5pt]
&\equiv 3^{2n}\sum_{k=0}^{n}{n\choose 3k}{n+k\choose k}{3k\choose 2k}{2k\choose k}\notag\\[5pt]
&\equiv 3^n c_n\pmod{8},\label{cong-proof-16}
\end{align}
where we have used \eqref{cong-proof-4}.

Furthermore, noting that for $k\ge 1$,
\begin{align*}
{4k+1\choose 2k}=\frac{2(4k+1)}{2k-1}{4k-1\choose 2k-2}\equiv 0\pmod{2},
\end{align*}
we have
\begin{align}
&3^{2n}\sum_{k=0}^{n}{2n\choose 3(2k+1)}{2n+2k+1\choose 2k+1}{3(2k+1)\choose 2(2k+1)}{2(2k+1)\choose 2k+1}\notag\\[5pt]
&=3^{2n}4n\sum_{k=0}^{n}\frac{1}{6k+3}{2n-1\choose 6k+2}{2n+2k+1\choose 2k+1}{6k+3\choose 4k+2}{4k+1\choose 2k}\notag\\[5pt]
&\equiv 4n(n-1)\notag\\[5pt]
&\equiv 0 \pmod{8}.\label{cong-proof-17}
\end{align}

Finally, combining \eqref{cong-proof-15}--\eqref{cong-proof-17} gives
\begin{align}
c_{2n}\equiv 3^nc_n\pmod{8}.\label{cong-proof-18}
\end{align}

By using the induction on $n$, \eqref{cong-proof-14} and \eqref{cong-proof-18}, we complete the proof of \eqref{cong-7}.
\qed

\vskip 5mm \noindent{\bf Acknowledgments.}
This work was supported by the National Natural Science Foundation of China (grant 12171370).

\end{document}